\documentclass[12pt,reqno]{article}
\usepackage{amsfonts}
\usepackage{amsfonts,mathrsfs}
\usepackage{CJK,fancyhdr,amscd}
\usepackage{anysize}
\usepackage{graphicx,epsfig}
\usepackage{amsthm,latexsym,amsmath,amssymb}
\usepackage[perpage,symbol*]{footmisc}
\usepackage{cite}
\usepackage[colorlinks=true]{hyperref}
\usepackage{color}
\newtheorem{thm}{Theorem}[section]

\newtheorem{lem}[thm]{Lemma}

\newcommand{\mc}{\mathcal}

\def\pf{\noindent{\bf Proof.} }

\makeatletter \@addtoreset{equation}{section} \makeatother
\def\beq{\begin{equation}}
\def\eeqp{\end{equation}}
\def\endproof{$\hfill\Box$}

\begin{document}
\baselineskip=20pt  \hoffset=-3cm \voffset=0cm \oddsidemargin=3.2cm
\evensidemargin=3.2cm \thispagestyle{empty}\vspace{10cm}
\title{\textbf{The existence of periodic solution for infinite dimensional  Hamiltonian systems}}
\author{\Large Weibing Deng$^{{\rm a}}$, $\quad$Wunming Han$^{{\rm a}}$,$\quad$ Qi Wang $^{{\rm a,b}}$}
\date{} \maketitle
\begin{center}
\it\scriptsize ${}^{\rm a}$ School of Mathematics and Statistics, Henan University, Kaifeng 475000, PR China\\
${}^{\rm b}$Department of Mathematics, Shandong University, Jinan, Shandong, 250100, PR China\\
\end{center}

\footnotetext[0]{$^b${\bf Corresponding author.} Supported by NNSF of China(11301148) and PSF of China(188576).}
\footnotetext[0]{\footnotesize\;{\bf E-mail address}: w$\_$deng8120@163.com(Weibing Deng),Wunming$\_$Han@sina.com (Wunming Han), Q.Wang@vip.henu.edu.cn. (Qi Wang).}

\date{} \maketitle

\noindent
{\bf Abstract:} {\small In this paper, we will consider a kind of infinite dimensional Hamiltonian system (HS), by the method of saddle point reduction, topology degree and the index defined in \cite{Wang-Liu-2015}, we will get the existence of periodic solution for (HS).}


\section{Introduction and main results}\label{Introduction and main results}
\subsection{Introduction of a kind of infinite dimensional Hamitonian system}
In this paper,  we will consider the following infinite dimensional Hamiltonian system
\[
\left\{\begin{array}{rl}
       \partial_t{u}-\Delta_xu&=H_v(t,x,u,v),\\
       -\partial_t{v}-\Delta_xv&=H_u(t,x,u,v),\\
       \end{array}
\right.\forall (t,x)\in\mathbb{R}\times\Omega, \eqno(HS)
\]
where $\Omega\subset\mathbb{R}^{N}$, $N\geq 1$ is a bounded domain with smooth boundary $\partial\Omega$ and $H:\mathbb{R}\times \overline{\Omega}\times\mathbb{R}^{2m}\to\mathbb{R}$ is a $C^1$ function,  $\partial_t:=\frac{\partial}{\partial t}$, $\Delta_x:=\sum_{i=1}^N\frac{\partial^2}{\partial x_i^2}$, $H_u:=\frac{\partial H}{\partial u}$ {  and  $H_v:=\frac{\partial H}{\partial v}$}. System like $(HS)$ are called unbounded Hamiltonian system, cf. Barbu\cite{Barbu-1995}, or infinite dimensional Hamiltonian system, cf. \cite{Bartsch-Ding-2002,Clement-Felemer-Mitidieri-1995,Clement-Felemer-Mitidieri-1997}. This systems arises in optimal control of systems governed by partial differential equations. See, e.g, Lions \cite{Lions-1971}, where the combination of the model $\partial_t-\Delta_x$ and its adjoint $-\partial_t-\Delta_x$ acts as a system for studying the control.  Br\'{e}zis and Nirenberg \cite{Brezis-Nirenberg-1978} considered  a special case of the system (HS):
\begin{equation}
\left\{\begin{array}{rl}
       \partial_t{u}-\Delta_xu&=-v^5+f,\\
       -\partial_t{v}-\Delta_xv&=u^3+g,\\
       \end{array}
\right.
\end{equation}
where $f,g\in L^\infty(\Omega)$, subject to the boundary condition $z(t,\cdot)|_{\partial\Omega}=0$ { on variable $x$} and the periodicity condition $z(0,\cdot)=z(T,\cdot)=0$ { on variable $t$} for a given $T>0$, where { $z=(u,v)$}. They obtained a solution $z$ with $u\in L^4$ and $v\in L^6$ by using Schauder's fixed point theorem. Cl\'{e}ment, Felemer and Mitidieri considered in \cite{Clement-Felemer-Mitidieri-1995} and \cite{Clement-Felemer-Mitidieri-1997}  the following system which is also a special case {  of (HS):}
\begin{equation}\label{Clement's problems}
\left\{\begin{array}{rl}
       \partial_t{u}-\Delta_xu&=|v|^{q-2}v,\\
       -\partial_t{v}-\Delta_xv&=|u|^{p-2}u,\\
       \end{array}
\right.
\end{equation}
with $\frac{N}{N+2}<\frac{1}{p}+\frac{1}{q}<1$.
Using their variational setting of Mountain Pass type, they proved that there is a $T_0>0$ such that, for each $T>T_0$, \eqref{Clement's problems} has at least one positive solution $z_T=(u_T,v_T)$ satisfying the boundary condition $z_T(t,\cdot)|_{\partial\Omega}=0$ for all $t\in (-T,T)$ and the periodicity condition $z_T(T,\cdot)=z_T(-T,\cdot)$ for all $x\in\overline{\Omega}$. Moreover, by passing to limit as $T\to \infty$ they obtained a positive homoclinic solution of \eqref{Clement's problems}. If the Hamiltonian function $H$ in (HS) can be displayed in the following form
\begin{equation}\label{a spectial form of Hamiltonian function}
H(t,x,u,v)=F(t,x,u,v)-V(x)uv,
\end{equation}
where $V\in C(\Omega,\mathbb{R})$, $H\in C^1(\mathbb{R}\times \overline{\Omega}\times\mathbb{R}^{2m},\mathbb{R})$, the system $(HS)$ will be rewritten as
\[
\left\{\begin{array}{rl}
       \partial_t{u}+(-\Delta_x+V(x))u&=F_v(t,x,u,v),\\
       -\partial_t{v}+(-\Delta_x+V(x))v&=F_u(t,x,u,v).\\
       \end{array}
\right.\eqno(HS.1)
\]
Bartsch and Ding\cite{Bartsch-Ding-2002} dealt with the system $(HS.1)$.
They established existence and multiplicity of { homoclinic} solutions of the type $z(t,x)\to 0 $ as $|t|+|x|\to\infty$ if $\Omega=\mathbb{R}^N$  and the type of $z(t,x)\to 0$ as $|t|\to\infty$ and $z(t,\cdot)|_{\partial\Omega}=0$ if $\Omega$ is bounded. Recently, there are several results on  system $(HS)$ and $(HS.1)$, cf.\cite{Ding-Lee-2005,Mao-Luan-2005,Wang-Xu-Zhang-Wang-2010,Wang-Liu-2015,Yang-Shen-Ding-2011,Zhang-Lv-Tang-2013,Zhang-Tang-Zhang-2015,Zhang-Tang-Zhang-2016}.

\subsection{Introduction of relative Morse index $(\mu_L(M),\upsilon_L(M))$}
In \cite{Wang-Liu-2015}, we developed the so called relative Morse index $(\mu_L(M),\upsilon_L(M))$ for (HS).
Let $I_m$ the identity map on $\mathbb R^m$ and
\begin{equation}
J=\left(\begin{matrix} 0 &-I_m\\I_m &0\end{matrix}\right),\;
N=\left(\begin{matrix} 0 &I_m\\I_m &0\end{matrix}\right),\;
\end{equation}

 \begin{equation}\label{operator L}
 L:=J\partial_t-N\Delta_x,
 \end{equation}
{ and denoted by $\nabla_z$ the gradient operator on variable} $z=(u,v)^T$,  then $(HS)$ with $T$-periodic and Dirichlet boundary conditions  can be rewritten as
 \[
 \left\{\begin{array}{rl}
                 L z&=\nabla_zH(t,x,z),\\
              z(t,x)&=z(t+T,x),\\
 z(t,\partial\Omega)&=0,
        \end{array}
\right.\forall (t,x)\in\mathbb{R}\times\Omega. \eqno(HS)
 \]
 Let  $\mathbf{H}:= L^2(S^1\times\Omega,\mathbb{R}^{2m})$, where $S^1=\mathbb{R}/T\mathbb{Z}$. Then $L$ is a self-adjoint operator acting in $\mathbf{H}$ with domain $D(L)$.
   The linearized system of the nonlinear system in (HS) at a solution $z=z(t,x)$ is the following system
  \begin{equation}\label{liu2} Ly=M(t,x)y\end{equation}
  with $M(t,x)=\nabla^2_zH(t,x,z(t,x))$ and  $\nabla^2_zH$  the Hessian of $H$ on the variable $z$.
  Denote by  $\mathbf{SM}^{2m}$ the set of all symmetric $2m\times 2m$ matrixes and  $\mathcal{M}:=C(S^1\times\bar{ \Omega},\mathbf{SM}^{2m})$. Denote by $\mathcal L_s(\mathbf H)$ the set of all bounded self-adjoint operators on $\mathbf H$. For any $M\in \mathcal{M}$, it is easy to see $M$ determines a  bounded self-adjoint operator  on $\mathbf{H}$, by
  \[
  z(t)\mapsto M(t,x)z(t), \;\forall z\in \mathbf{H},
\]
we still denote this operator by $M$. Thus, we have $\mathcal M\subset\mathcal L_s(\mathbf H)$. In \cite{Wang-Liu-2015}, for any $B\in\mathcal{L}_s(\mathbf H)$, we  defined  the relative Morse index pair
\begin{equation}
(\mu_L(B), \upsilon_L(B))\in\mathbb{Z}\times\mathbb{Z}^*,
\end{equation}
where $\mathbb{Z}$ and $\mathbb{Z}^*$ denote the set of all integers and  non-negative integers respectively. Then, we got the relationship between the index $\mu_L(B)$ and other indexes. Spectrally, with the relationship  between the index $\mu_L(B)$ and spectral flow, we have the following equality which  will be used in this paper. For  $B_1,B_2\in\mathcal L_s(\mathbf H)$, $B_1\le B_2$ means that $B_2-B_1$ is semi-positive definite. Then, if $B_1\le B_2$, we have
\begin{equation}\label{eq-indexaddition}
\mu_L(B_2)-\mu_L(B_1)=\sum_{s\in[0,1)}\nu_L(sB_2+(1-s)B_1).
\end{equation}
By assuming some twisted conditions of the  asymptotically linear Hamiltonian function, we studied the existence and multiplicity of $(HS)$ in \cite{Wang-Liu-2015}.
\subsection{Main results}
In this paper, we don't need the Hamiltonian function $H$ to be $C^2$ continuous and without assuming $H$ satisfying the twisted conditions, by the method of topology degree, saddle point reduction and the index $(i_L(B),\nu_L(B))$ defined in \cite{Wang-Liu-2015}, we have the following results.
\begin{thm}\label{thm-1}
Assume   $H$ satisfies  the following conditions.\\
($H_1$) $H\in C^1(S^1\times\bar\Omega\times\mathbb R^{2m},\mathbb R)$ and there exists $l_H>0$, such that
\[
|H^{\prime}_z(t,x,z+y)-H^{\prime}_z(t,x,z)|\leq l_H|y|,\;\;\forall (t,x)\in S^1\times\bar\Omega,\;z, y\in\mathbb{R}^{2m}.
\]
($H^\pm_2$) There exists $M_1,\;M_2,\;K>0$, $B\in \mc M$, such that
\[
H^\prime(t,x,z)=B(t,x)z+r(t,x,z),
\]
with
\[
|r(t,x,z)| \leq M_1,\;\;\forall (t,x,z)\in S^1\times\bar\Omega\times\mathbb{R}^{2m},
\]
and
\begin{equation}\label{eq-condition of r in ab-thm 3}
\pm(r(t,x,z),z)_{\mathbb{R}^{2m}}\geq M_2|z|_{\mathbb{R}^{2m}},\;\forall (t,x)\in S^1\times\bar\Omega,\;\|z\|_{\mathbb{R}^{2m}}>K.
\end{equation}
 Then (HS) has at least one solution.
\end{thm}
In Theorem \ref{thm-1}, we don't need $B$ to be  non-degenerate. If we assume some  non-degenerate property of $B$, the rest item $r$ can be relaxed and we have the following result.
\begin{thm}\label{thm-2}
Assume $H$ satisfying condition ($H_1$) and the following condition\\
($H_3$) There exists $B\in C(S^1\times\bar\Omega\times\mathbb{R}^{2m},\mathbf{SM}^{2m})$  such that\\
\[
H^\prime(t,z)=B(t,x,z)z+r(t,x,z),\;\forall (t,x,z)\in S^1\times\bar\Omega\times\mathbb{R}^{2m},
\]
with
\[
r(t,x,z)=o(z),\;{\rm uniformly \;for }|z|\to \infty.
\]
($H_4$) There exist  $B_1,\;B_2\in \mathcal M$
satisfying
\[
i_L(B_1)=i_L(B_2),\;\;\nu_L(B_2)=0,
\]
and
\[
B_1(t,x)\leq B(t,x,z)\leq B_2(t,x),\;\forall(t,x,z)\in S^1\times\bar\Omega\times\mathbb{R}^{2m}.
\]
Then (HS) has at least one solution.
\end{thm}
\section{Preliminarys and the proof of our main results}
Before the proof of Theorem \ref{thm-1} and Theorem \ref{thm-2}, we need some preliminarys. Firstly, we need the following Lemma.
\begin{lem}\label{Lemma-The spectral of L}\cite[Lemma 2.1]{Wang-Liu-2015}
For simplicity, let $T=2\pi$ and $\sigma(-\Delta)=\{\mu_l\}_{l\geq 1}$, we have
\[
\sigma(L)=\sigma_p(L)=\{\pm(k^2+\mu^2_l)^{1/2}\}_{k\in\mathbb{Z},l\in\mathbb{N}},
\]
where { $\sigma_p(L)$ denotes the eigenvalue set of $L$ on $\mathbf{H}$}. That is to say $L$ has only eigenvalues. More over every eigenvalue in $\sigma_p(L)$ has $2m$ dimensional eigenspace.
\end{lem}
Secondly, since $H$ satisfies condition ($H_1$), the map
\[
z\mapsto\int_{S^1\times\bar\Omega}H(t,x,z(t,x))dtdx,\;\forall z\in \mathbf H,
\]
define a functional on $\mathbf H$, without confusion, we still denote it by $H$.
It is easy to see $H\in C^1(\mathbf H,\mathbb{R})$, with
\[
(H^\prime(z),y)_{\mathbf H}=\int_{S^1\times\bar\Omega}(H^\prime_z(t,x,z(t,x)),y(t))dtdx,\;\forall z,y\in\mathbf H,
\]
 and $H^\prime$ is  Lipschitz continuous with
\begin{equation}\label{eq-lip of Phi}
\|H^\prime(z+y)-H^\prime(z)\|_\mathbf H\leq l_H\|y\|_\mathbf H,\;\forall,y,z\in\mathbf H.
\end{equation}
Thus (HS) can be regard as an operator equation on $\mathbf H$.

\noindent{\bf Proof of Theorem \ref{thm-1}.}  Now, we consider the case of ($H^-_2$). From Lemma \ref{Lemma-The spectral of L}, $L$ has compact resolvent, since $B\in \mathcal{M}\subset \mathcal L_s{\mathbf (H)}$, $0$ is at most an isolate point spectrum with finite dimensional eigenspace, that is to say there exists $\varepsilon_0>0$ and small enough, such that $(-\varepsilon_0,0)\cap \sigma(L-B)=\emptyset$.   For any $\varepsilon\in(0,\varepsilon_0)$ and $\lambda\in[0,1]$, consider the following two-parameters equation
\[
(\varepsilon\cdot I+L-B)z=\lambda r(t,x,z),\eqno{(HS_{\varepsilon,\lambda})}
\]
with $I$ the identity map on $\mathbf H.$
If $\varepsilon=0$ and $\lambda=1$, it is (HS). We divide the following proof into four steps.\\
{\bf Step 1. There exists a constant $C$ independent of $\varepsilon$ and $\lambda$, such that if $z_{\varepsilon,\lambda}$ is a solution of ($HS_{\varepsilon,\lambda}$),
\[
\varepsilon\|z_{\varepsilon,\lambda}\|_\mathbf H\leq C,\;\forall (\varepsilon,\lambda)\in(0,\frac{\varepsilon_0}{2})\times[0,1].
\]
}
Since $(-\varepsilon_0,0)\cap \sigma(L-B)=\emptyset$, we have $(\varepsilon-\varepsilon_0,\varepsilon)\cap \sigma(\varepsilon\cdot I+L-B)=\emptyset$. Consider the orthogonal splitting
\[
\mathbf H=\mathbf H^-_{\varepsilon\cdot I+L-B}\oplus\mathbf H^-_{\varepsilon\cdot I+L-B},
\]
where $\varepsilon\cdot I+L-B$ is negative definite on $\mathbf H^-_{\varepsilon\cdot I+L-B}$, and positive define on $\mathbf H^+_{\varepsilon\cdot I+L-B}$. Thus , if $z\in\mathbf H$, we have the splitting
\[
z=x+y,
\]
with $x\in\mathbf H^-_{\varepsilon\cdot I+L-B}$ and $y\in\mathbf H^+_{\varepsilon\cdot I+L-B}$.
If $z_{\varepsilon,\lambda}$ is a solution of ($HS_{\varepsilon,\lambda}$) with its splitting $z_{\varepsilon,\lambda}=x_{\varepsilon,\lambda}+y_{\varepsilon,\lambda}$ defined above, then we have
\[
((\varepsilon\cdot I+L-B)z_{\varepsilon,\lambda},y_{\varepsilon,\lambda}-x_{\varepsilon,\lambda})_\mathbf H=\lambda(r(t,z_{\varepsilon,\lambda}),y_{\varepsilon,\lambda}-x_{\varepsilon,\lambda})_\mathbf H.
\]
Since $(\varepsilon-\varepsilon_0,\varepsilon)\cap \sigma(\varepsilon\cdot I+L-B)=\emptyset$, we have
\[
((\varepsilon\cdot I+L-B)z_{\varepsilon,\lambda},y_{\varepsilon,\lambda}-x_{\varepsilon,\lambda})_\mathbf H\geq \min\{\varepsilon_0-\varepsilon,\varepsilon\}\|z_{\varepsilon,\lambda}\|^2_\mathbf H.
\]
Since $r$ is bounded, for $(\varepsilon,\lambda)\in(0,\frac{\varepsilon_0}{2})\times[0,1]$, we have
\[
C\|z_{\varepsilon,\lambda}\|_\mathbf H \geq\lambda(r(t,x,z_{\varepsilon,\lambda}),y_{\varepsilon,\lambda}-x_{\varepsilon,\lambda})_\mathbf H\geq \varepsilon \|z_{\varepsilon,\lambda}\|^2_\mathbf H.
\]
Therefor, we have
\[
\varepsilon\|z_{\varepsilon,\lambda}\|_\mathbf H\leq C,\;\forall (\varepsilon,\lambda)\in(0,\frac{\varepsilon_0}{2})\times[0,1].
\]
 {\bf Step 2. For any $(\varepsilon,\lambda)\in(0,\frac{\varepsilon_0}{2})\times[0,1]$, ($HS_{\varepsilon,\lambda}$) has at least one solution.}
Here, we use the topology degree theory. Since $0\notin\sigma(\varepsilon\cdot I+L-B)$, ($HS_{\varepsilon,\lambda}$) can be rewritten as
\[
z=\lambda(\varepsilon\cdot I+L-B)^{-1}r(t,x,z).
\]
Denote by $f(\varepsilon,\lambda,z):=\lambda(\varepsilon\cdot I+L-B)^{-1}r(t,x,z)$ for simplicity.
From the compactness of $(\varepsilon\cdot I+L-B)^{-1}$ and condition ($H^-_2$),  Leray Schauder degree theory can be used to the map
\[
z\mapsto z-f(\varepsilon,\lambda,z),\;z\in\mathbf{H}.
\]
From the result received in Step 1, we have
\begin{align*}
deg(I-f(\varepsilon,\lambda,\cdot),B(R(\varepsilon),0),0)&\equiv deg(I-f(\varepsilon,0,\cdot),B(R(\varepsilon),0),0)\\
                                                                                   &=deg(I,B(R(\varepsilon),0),0)\\
                                                                                   &=1,
\end{align*}
where $R(\varepsilon)>\frac{C}{\varepsilon}$ is a constant only depends on $\varepsilon$,
and $B(R(\varepsilon),0):=\{z\in\mathbf H |\|z\|_\mathbf H<R(\varepsilon)\}$.\\
{\bf Step 3. For $\lambda=1$, $\varepsilon\in(0,\varepsilon_0/2)$,  denote by  $z_\varepsilon$ one of the solutions of ($HS_{\varepsilon,1}$). We have  $\|z_\varepsilon\|_\mathbf H\leq C$. In this step,  $C$  denotes various constants independent of $\varepsilon$.} \\
From the boundedness received in Step 1, we have
\begin{equation}\label{eq-boundedness of z-1}
\|(L-B)z_\varepsilon\|_\mathbf H=\|\varepsilon z_\varepsilon-r(t,x,z_\varepsilon)\|_\mathbf H\leq C.
\end{equation}
Now, consider the orthogonal splitting
\[
\mathbf L=\mathbf H^0_{L-B}\oplus\mathbf H^\bot_{L-B},
\]
where $L-B$ is zero definite on $\mathbf H^0_{L-B}$,  and $\mathbf H^\bot_{L-B}$ is the orthonormal complement space of $\mathbf H^0_{L-B}$. Let $z_\varepsilon=u_\varepsilon+v_\varepsilon$ with $u_\varepsilon\in \mathbf H^0_{L-B}$ and $v_\varepsilon\in \mathbf H^\bot_{L-B}$. Since $0$ is an isolated point in $\sigma(L-B)$, from \eqref{eq-boundedness of z-1}, we have
\begin{equation}\label{eq-boundedness of z-2}
\|v_\varepsilon\|_\mathbf H\leq C
\end{equation}
Additionally, since $r(t,x,z_\varepsilon)$ and $v_\varepsilon$ are bounded in $\mathbf H$, we have
\begin{align}\label{eq-boundedness of z-3}
(r(t,x,z_\varepsilon),z_\varepsilon)_\mathbf H&=(r(t,x,z_\varepsilon),v_\varepsilon)_\mathbf H+(r(t,x,z_\varepsilon),u_\varepsilon)_\mathbf H\nonumber\\
                                           &=(r(t,x,z_\varepsilon),v_\varepsilon)_\mathbf H+(\varepsilon z_\varepsilon+(L-B)z_\varepsilon,u_\varepsilon)_\mathbf H\nonumber\\
                                                                                      &=(r(t,x,z_\varepsilon),v_\varepsilon)_\mathbf H+\varepsilon(u_\varepsilon,u_\varepsilon)_\mathbf H\nonumber\\
                                           &\geq C.
\end{align}
On the other hand, from \eqref{eq-condition of r in ab-thm 3} in ($H^-_2$),  we have
\begin{align}\label{eq-boundedness of z-4}
(r(t,x,z_\varepsilon),z_\varepsilon)_\mathbf H&=\int_{S^1\times\bar\Omega}(r(t,x,z_\varepsilon),z_\varepsilon)dtdx+\int_{S^1\times\bar\Omega/S^1\times\bar\Omega(K))}(r(t,x,z_\varepsilon),z_\varepsilon)dtdx\nonumber\\
&\leq-M_2\int_{S^1\times\bar\Omega(K)}|z_\varepsilon|^2dtdx+C\nonumber\\
&\leq-M_2\|z_\varepsilon\|^2_\mathbf H+C,
\end{align}
where $S^1\times\bar\Omega(K):=\{(t,x)\in S^1\times\bar\Omega||z_\varepsilon|>K\}$.
 From \eqref{eq-boundedness of z-3} and \eqref{eq-boundedness of z-4}, we have proved the boundedness of $\|z_\varepsilon\|_\mathbf H$.\\
{\bf Step 4. Passing to a sequence of $\varepsilon_n\to 0$,  there exists $z\in\mathbf H$
such that
\[
\displaystyle\lim_{\varepsilon_n\to 0}\|z_{\varepsilon_n}-z\|_\mathbf H=0.
\]
}
Here, we will use the method of saddle point reduction. Since  $\sigma(L)$  has only isolate finite dimensional eigenvalues and from condition ($H_1$),  we can assume $\pm l_H\notin\sigma(L)$. That is to say there exists $\delta>0$ such that
\[
(-l_H-\delta,-l_H+\delta)\cap\sigma(L)=(l_H-\delta,l_H+\delta)\cap\sigma(L)=\emptyset.
\]
Denote $E_L$ the spectrum measure of $L$ and definite the projections on $\mathbf H$ by
\begin{equation}\label{eq-projections}
P^0_{L,l_H}:=\int^{l_H}_{-l_H}dE_L(z),\;P^{\bot}_{L,l_H}:=I-P^0_{L,l_H},
\end{equation}
where $I$ is the identity map on $\mathbf H$.
Correspondingly, consider the splitting of $\mathbf H$ by
\begin{equation}\label{eq-splitting of L}
\mathbf H=\mathbf H^0_{L,l_H}\oplus\mathbf H^\bot_{L,l_H}
\end{equation}
with $\mathbf H^\star_L:=P^\star_{L,l_H}\mathbf H$ ($\star=0,\bot$). Without confusion, we rewrite $P^\star:=P^\star_{L,l_H}$ and $\mathbf H^\star:=\mathbf H^\star_{L,l_H}$ for simplicity ($\star=0,\bot$). Denote $L^\star=L|_{\mathbf H^\star}$ and $z^\star=P^\star z$, forall $z\in\mathbf H$. thus we have $L^\bot$ has bounded inverse on $\mathbf H^\bot$ and
\[
\|(L^\bot)^{-1}\|\leq\frac{1}{l_H+\delta}.
\]
Let $\varepsilon^\prime:=\min\{\varepsilon_0,\delta\}$, for $\varepsilon\in (0,\frac{\varepsilon^\prime}{2})$, denote by $L_\varepsilon:=\varepsilon+L$. Then $L_\varepsilon$ has the same invariant subspace with $L$, so we can also denote by $L^\star_\varepsilon:=L_\varepsilon|_{\mathbf H^\star}$ ($\star=0,\bot$), and we have
\begin{equation}\label{eq-bound of Apm}
\|(L^\bot_\varepsilon)^{-1}\|\leq\frac{1}{l_H+\delta/2}.
\end{equation}
Since $z_\varepsilon$ satisfies ($HS_{\varepsilon,1}$), so we have
\[
L^\bot_\varepsilon z^\bot_\varepsilon =P^\bot H^\prime(z^\bot_\varepsilon+z^0_\varepsilon),
\]
and
\begin{equation}\label{eq-splitting of OE}
z^\bot_\varepsilon =(L^\bot_\varepsilon)^{-1}P^\bot\Phi^\prime(z^\bot_\varepsilon+z^0_\varepsilon).
\end{equation}
Since $\mathbf H^0$  is a finite dimensional space and $\|z_\varepsilon\|_\mathbf H\leq C$,  there exists a sequence $\varepsilon_n\to 0$ and $z^0\in \mathbf H^0$, such that
\[
\displaystyle\lim_{n\to\infty}z^0_{\varepsilon_n}=z^0.
\]
 For simplicity, we rewrite $z^\star_n:=z^\star_{\varepsilon_n}$($\star=\bot,0$), $L_n:=\varepsilon_n+L$ and $L^\bot_n:=L^\bot_{\varepsilon_n} $.
 So, we have
 \begin{align*}
\|z^\bot_n-z^\bot_m\|_\mathbf H=&\|(L_n^\bot)^{-1}P^\bot \Phi^\prime(z_n)-(L_m^\bot)^{-1}P^\bot \Phi^\prime(z_m)\|_\mathbf H \\
\leq&\|(L_n^\bot)^{-1}P^\bot (\Phi^\prime(z_n)-\Phi^\prime(z_m))\|_\mathbf H+\|((L_n^\bot)^{-1}-(L_m^\bot)^{-1})P^\bot \Phi^\prime(z_m)\|_\mathbf H\\
\leq&\frac{l_H}{l_H+\delta/2}\|z_n-z_m\|_\mathbf H+\|((L_n^\bot)^{-1}-(L_m^\bot)^{-1})P^\bot \Phi^\prime(z_m)\|_\mathbf H\\
\leq&\frac{l_H}{l_H+\delta/2}(\|z^\bot_n-z^\bot_m\|_\mathbf H+\|z^0_n-z^0_m\|_\mathbf H)+\|((L_n^\bot)^{-1}-(L_m^\bot)^{-1})P^\bot \Phi^\prime(z_m)\|_\mathbf H.
\end{align*}
Since $(L_n^\bot)^{-1}-(L_m^\bot)^{-1}=(\varepsilon_m-\varepsilon_n)(L_n^\bot)^{-1}(L_m^\bot)^{-1}$ and $z_n$ are bounded in $\mathbf H$, we have
\[
\|((L_n^\bot)^{-1}-(L_m^\bot)^{-1})P^\bot \Phi^\prime(z_m)\|_\mathbf H=o(1),\;\;n,m\to\infty.
\]
So we have
\[
\|z^\bot_n-z^\bot_m\|_\mathbf H\leq \frac{2l_H}{\delta}\|z^0_n-z^0_m\|_\mathbf H+o(1),\;\;n,m\to\infty,
\]
therefor, there exists $z^\bot\in\mathbf H^\bot$, such that $\displaystyle\lim_{n\to\infty}\|z^\bot_n- z^\bot\|_\mathbf H=0$. Thus, we have
\[
\displaystyle\lim_{n\to\infty}\|z_{\varepsilon_n}-z\|_\mathbf H=0,
\]
with $z=z^\bot+z^0$.
Last, let $n\to\infty$ in ($HS_{\varepsilon_n,1}$), we have $z$ is a solution of ($HS$).
\endproof

Before the proof of Theorem \ref{thm-2}, we need the following Lemma.
\begin{lem}\label{lem-0 has a positive distance from sigma(A-B)}
Let $B_1,B_2\in \mathcal{L}_s(\mathbf H)$ with $B_1\leq B_2,\;\mu_L(B_1)=\mu_L(B_2),\; {\rm and}\; \nu_L(B_2)=0$, then there exists $\varepsilon>0$, such that for all $B\in\mathcal{L}_s(\mathbf H)$ with
\[
 B_1\leq B \leq B_2,
\]
we have
\[
\sigma(L-B)\cap (-\varepsilon,\varepsilon)=\emptyset.
\]
\end{lem}
{\noindent}{\bf Proof.} For the property of $\mu_L(B)$, we have $\nu_L(B_1)=0$. So there is $\varepsilon>0$, such that
\[
\mu_L(B_{1,\varepsilon})=\mu_L(B_1)=\mu_L(B_2)=\mu_L(B_{2,\varepsilon}),
\]
with $B_{*,\varepsilon}=B_*+\varepsilon\cdot I,(*=1,2)$. Since  $B_{1,\varepsilon}\leq B-\varepsilon I<B+\varepsilon I\leq B_2'$. It follows that $\mu_L(B-\varepsilon I)=\mu_L(B+\varepsilon I)$. Note that by \eqref{eq-indexaddition}
\begin{align*}
\sum_{-\varepsilon < t \le  \varepsilon } \nu_L(B-t \cdot I)&=\mu_L(B+\varepsilon)-\mu_L(B-\varepsilon)\\
&=0.
\end{align*}
We have $0\notin \sigma(L-B-\eta),\;\forall \eta\in(-\varepsilon,\varepsilon)$, thus the proof is complete.\endproof

 \noindent{\bf Proof of Theorem \ref{thm-2}.}  Consider the following one-parameter equation
 \[
 Lz=(1-\lambda) B_1z+\lambda H^\prime(z),\eqno(HS_\lambda)
 \]
 with  $\lambda\in[0,1]$.
 Denote by
\[
\Phi_\lambda(z)=\frac{1-\lambda}{2}(B_1z,z)_{\mathbf L}+\lambda \Phi(z),\;\forall z\in \mathbf L.
 \]
 Since $H$ satisfies condition ($H_1$) and $B_1\in C(S^1\times\bar\Omega,\mathbf{SM}^{2m})$, we have $\Phi¡¯_\lambda:\mathbf H\to\mathbf H$ is Lipschitz continuous , and there exists   $l^\prime>0$ independed of $\lambda$ such that $l^\prime\notin\sigma(L)$ and
 \[
 \|\Phi^\prime_\lambda(z+h)-\Phi^\prime_\lambda(z)\|_{\mathbf H}\leq l^\prime\|h\|_{\mathbf H},\;\forall z,h\in\mathbf L,\lambda\in[0,1].
 \]
  Now, replace $l_H$ by $l^\prime$ in \eqref{eq-projections}, we have  the projections $P^\star_{L,l^\prime}$($\star=\bot,0$) and the splitting
  \[
  \mathbf H=\mathbf H^\bot_{L,l^\prime}\oplus\mathbf H^0_{L,l^\prime},
   \]
   with $\mathbf H^\star_{L,l^\prime}=P^\star_{L,l^\prime}\mathbf H$($\star=\bot,0$). Thus $L^\bot$ has bounded inverse on $\mathbf H^\bot_{L,l^\prime}$ with
   \[
   \|(L^\bot)^{-1}\|<\frac{1}{l^\prime+c},
   \]
   for some $c>0$. Without confusion, we still use $z^\bot$  and $z^0$ to represent the splitting
 \[
 z=z^\bot+z^0,
 \]
 with $z^\star\in\mathbf H^\star_{L,l^\prime}$($*=\bot,0$). Now, we derive the following proof into three steps and $C$  denotes various constants independent of $\lambda$. \\
 {\bf Step 1. If $z$ is a solution of ($HS_\lambda$), then we have $\|z^\bot(z^0)\|_\mathbf H\leq C\|z^0\|_{\mathbf H}+C$}\\
 Since $L z=\Phi_\lambda^\prime(z)$, we have
 \[
 z^\bot=(L^\bot)^{-1}P^\bot_{L,l^\prime}\Phi^\prime_\lambda(z)
  \]
 \begin{align*}
\|z^\pm(x)\|_{\mathbf H}&=\|(L^\pm)^{-1}P^\pm_L \Phi^\prime_\lambda(z^\bot(z^0)+z^0)\|_{\mathbf H}\\
                &\leq\frac{1}{l^\prime+c}\|\Phi^\prime_\lambda(z^\bot(z^0)+z^0)\|_{\mathbf H}\\
                &\leq\frac{1}{l^\prime+c}\|\Phi^\prime_\lambda(z^\bot(z^0)+z^0)-\Phi^\prime_\lambda(0)\|_{\mathbf H}+\frac{1}{l^\prime+c}\|\Phi^\prime_\lambda(0)\|_{\mathbf H}\\
                &\leq\frac{l^\prime}{l^\prime+c}(\|z^\bot(z^0)\|_{\mathbf H}+\|z^0\|_\mathbf H)+\frac{1}{l^\prime+c}\|\prime_\lambda(0)\|_{\mathbf H}.
\end{align*}
So we have $\|z^\bot(z^0)\|_\mathbf H\leq \frac{l^\prime}{c}\|z^0\|_{\mathbf H}+\frac{1}{c}\|\Phi^\prime_\lambda(0)\|_{\mathbf H}$.
Thus, we have prove this step.

 \noindent{\bf Step 2. We claim that the set of all the solutions ($z,\lambda$) of ($HS_\lambda$) are a priori bounded. }\\
 If not, assume there exist $\{(z_n,\lambda_n)\}$ satisfying ($HS_\lambda$) with $\|z_n\|_{\mathbf H}\to\infty$. Without lose of generality, assume $\lambda_n\to\lambda_0\in[0,1]$.
  From step 1,  we have $\|z^0_n\|_\mathbf L\to\infty$. Denote by
 \[
 y_n=\frac{z_n}{\|z_n\|_\mathbf H},
 \]
 and $\bar{B}_n:=(1-\lambda_n)B_1+\lambda_nB(t,z_n)$, we have
\begin{equation}\label{eq-the equation of yn}
 Ay_n=\bar{B}_ny_n+\frac{o(\|z_n\|_\mathbf H)}{\|z_n\|_\mathbf H}.
 \end{equation}
Decompose $y_n=y^{\bot}_n+y^0_n$ with $y^\star_n=z^\star_n/\|z_n\|_\mathbf H$, we have
\begin{align*}
\|y^0_n\|_\mathbf H&=\frac{\|z^0_n\|_\mathbf H}{\|z_n\|_\mathbf H}\\
                   &\geq\frac{\|z^0_n\|_\mathbf H}{\|z^0_n\|_\mathbf H+\|z^\bot\|_\mathbf H}\\
                   &\geq\frac{\|z^0\|_\mathbf H}{C\|z^0\|_\mathbf H+C}.
\end{align*}
That is to say
 \begin{equation}\label{eq-y0n not to 0}
\|y^0_n\|_\mathbf H\geq C>0,
 \end{equation}
 for  $n$ large enough.
 Since $B_1(t)\leq B(t,z)\leq B_2(t)$, we have $B_1\leq \bar{B}_n\leq B_2$.
 Let $\mathbf H=\mathbf H^+_{L-\bar B_n}\bigoplus \mathbf H^-_{L-\bar B_n}$ with $L-\bar B_n$ is positive and negative define on $\mathbf H^+_{L-\bar B_n}$ and $\mathbf H^-_{L-\bar B_n}$ respectively. Re-decompose $y_n=\bar{y}^+_n+\bar{y}^-_n$ respect to $\mathbf H^+_{L-\bar B_n}$ and $\mathbf H^-_{L-\bar B_n}$. From ($H_4$) and \eqref{eq-the equation of yn}, we have
 \begin{align}\label{eq-y0n to 0}
 \| y^0_n\|^2_\mathbf H&\leq  \| y_n\|^2_\mathbf H\nonumber\\
                     &\leq C ((A-\bar{B}_n)y_n,\bar{y}^+_n+\bar{y}^-_n)_\mathbf H\nonumber\\
                     &\leq  \frac{o(\|z_n\|_\mathbf H)}{\|z_n\|_\mathbf H}\|y_n\|_\mathbf H.
 \end{align}
 Since $\|z_n\|_\mathbf H\to \infty$ and $\|y_n\|_\mathbf H=1$, we have $\|y^0_n\|_\mathbf H\to 0$ which  contradicts to \eqref{eq-y0n not to 0}, so we have $\{z_n\}$ is bounded.

 \noindent{\bf Step 3. By Leray-Schauder degree,  there is a solution of (HS)}.\\
    Since the solutions of ($HS_\lambda$ )are bounded, there is a number $R>0$ large eoungh,
 such that all of the solutions $z_\lambda$ of ($HS_\lambda$ ) are in the ball $B(0,R):=\{z\in \mathbf L| \|z\|_\mathbf H<R\}$.  So we have the Larey-Schauder degree
 \[
deg (I-(L-B_1)^{-1}(\Phi^\prime(z)-B_1z),B(0,R)\cap,0)= deg (I,B(0,R),0)=1.
 \]
  That is to say (HS) has at least one solution.
  \section{Further results}
In the system of (HS), Lemma \ref{Lemma-The spectral of L} played an important role to keep the Leray-Schauder degree valid,  if we change the Dirichlet boundary condition $z(t,\partial\Omega)=0$ in (HS) to Neumann boundary condition $\displaystyle\frac{\partial z}{\partial n }(t,\partial\Omega)=0$, we will also have Lemma \ref{Lemma-The spectral of L}, thus Theorem \ref{thm-1} and Theorem \ref{thm-2}  will also be true for  Neumann boundary condition.

What we want to say in this section is $\Omega=\mathbb{R}^N$. Generally, the operator $-\Delta_x$ on $L^2(\mathbb{R}^N,\mathbb{R})$ doesn't have compact inversion,  then the results in Lemma \ref{Lemma-The spectral of L} will not be true. Thus our Maslov type index theory defined in \cite{Wang-Liu-2015} will not work. But if the Hamiltonian function $H$ can be displayed in the following form
\begin{equation}
H(t,x,u,v)=F(t,x,u,v)-V(x)uv,
\end{equation}
then system ($HS$) will be rewritten as systems
\[
\left\{\begin{array}{rl}
       \partial_t{u}+(-\Delta_x+V(x))u&=F_v(t,x,u,v),\\
       -\partial_t{v}+(-\Delta_x+V(x))v&=F_u(t,x,u,v).\\
       \end{array}
\right.\eqno(HS.1)
\]
We have the following result.
\begin{lem} \cite[Lemma 6.10]{Ding-2007}.
If  the function $V(x)$ satisfies the following conditions:\\
($V_1$) $V\in C(\mathbb{R}^N,\mathbb{R})$ and $\inf_{x\in\mathbb{R}^N}V(x)>0$.\\
($V_2$) There exists $l_0>0$ and $M>0$ such that
\[
\displaystyle\lim_{|y|\to\infty}meas(\{x\in\mathbb{R}^N:|x-y|\leq l_0,\;V(x)\leq M\})=0,
\]
where $meas(\cdot)$ denotes the Lebesgue measure in $\mathbb{R}^N$. Then we have
\[
\sigma_{e}(-\Delta_x+V(x))\subset [M,+\infty),
\]
where $\sigma_e(A)$ denotes the essential spectrum of operator $A$.
\end{lem}  If we redefine the operator $L$ as
\[
L:=J\partial_t-N(\Delta_x-V(x)).
\]
We have the following result.
\begin{lem} If $\sigma_{e}(-\Delta_x+V(x))\subset [M,+\infty)$, then
\[
\sigma_{e}(L)\cap (-M,M)=\emptyset.
\]
\end{lem}
\pf Let $E(z)$ be the spectrum measure of $-\Delta_x+V(x)$,  for any $\delta>0$ small enough, define the  following projection on $L^2(\mathbb{R}^N,\mathbb R)$,
\[
P_{-\Delta_x+V(x),\delta}=\int^{M-\delta}_0dE(z).
\]
We have the following orthogonal splitting
\[
L^2(\mathbb{R}^N,\mathbb R)=L^2(\delta)\oplus (L^{2,\bot}(\delta)),
\]
where $L^2(\delta):=P_{-\Delta_x+V(x),\delta} L^2(\mathbb{R}^N,\mathbb R)$ and $(L^{2,\bot}(\delta))$ is its orthogonal complement. So we have
\[
-\Delta_x+V(x)|_{L^2(\delta)}\leq M-\delta,\;\;-\Delta_x+V(x)|_{(L^2(\delta))^\bot}\geq M-\delta.
\]
Since
\begin{align*}
L^2(S^1\times\mathbb R^N,\mathbb R^{2m})&=L^2(S^1,L^2(\mathbb R^N,\mathbb R^{2m}))\\
                                                                       &=L^2(S^1,L^2(\mathbb R^N,\mathbb R)^{2m})\\
                                                                       &=L^2(S^1,(L^2(\delta)\oplus L^{2,\bot}(\delta))^{2m})\\
                                                                       &=L^2(S^1,L^2(\delta)^{2m}\oplus L^{2,\bot}(\delta)^{2m})\\
                                                                       &=L^2(S^1,L^2(\delta)^{2m})\oplus L^2(S^1,L^{2,\bot}(\delta)^{2m}),
\end{align*}
 $L^2(S^1,L^2(\delta)^{2m})$ and $L^2(S^1,L^{2,\bot}(\delta)^{2m})$ are invariant subspaces of $L$, let
 \[
 L_1:=L|_{L^2(S^1,L^2(\delta)^{2m})},\;L_2:=L|_{L^2(S^1,L^{2,\bot}(\delta)^{2m})}.
 \]
 Corresponding to the splitting of $L^2(S^1\times\mathbb R^N,\mathbb R^{2m})$, we have
\[
L=\left(\begin{matrix} L_1 &0\\0 & L_2\end{matrix}\right).
\]
So we have
\[
\sigma(L)=\sigma(L_1)\cup \sigma(L_2).
\]
With the similarly method in Lemma \ref{Lemma-The spectral of L}, we can prove  $\sigma_e(L_1)=\emptyset$,  so
\[
\sigma_e(L)=\sigma_e(L_2).
\]
Now, we will prove $\sigma(L_2)\cap(-M,M)=\emptyset$. For any $\lambda\in\sigma(L_2)$, we have $z_n\in L_2$ with $\|z_n\|=1$, such that
\[
\|L_2z_n-\lambda z_n\|\to 0.
\]
So we have
\[
(Lz_n-\lambda z_n,Nz_n)\to 0,
\]
that is to say
\[
(-J\frac{\partial}{\partial t}z_n-N(\Delta-V(x))z_n,Nz_n)-\lambda(z_n,Nz_n)\to 0.
\]
Since $(-J\frac{\partial}{\partial t}z_n,Nz_n)\equiv 0$ and from $z_n\in L_2$ we have $(-N(-\Delta-V(x))z_n,Nz_n)\geq M$,
so we have $|\lambda|\geq M$. Thus, we have finished the proof.\endproof

 With this lemma and the index defined in \cite{Wang-Liu-2016}, we can also get the similar results as Theorem \ref{thm-1} and \ref{thm-2}, we omit them here.

 {\bf Acknowledgements.} The author of this paper sincerely thanks the referee for his/her
careful reading and valuable comments and suggestions on the first manuscript of this
paper.
 
\end{document}